# The Gion Shrine Problem:
# A Solution in Geometry

Melissa Holly

# Table of Contents

# List of Figures and Tables

# Abstract

Eighteenth century Japan was a time of isolation and peace, where education and the arts blossomed. Originally posted before 1749 by an unknown author, the *sangaku* (mathematical tablet) that became known as the Gion Shrine problem, has puzzled mathematicians for over two centuries.  Although solutions have been suggested, some of them require mathematics not known in Japan at the time the Gion Shrine problem was written.  The problem is still considered unsolved.  Using the geometry of Edo Period Japan, this paper's solution is immersed in the culture of that timeframe. As identified in this paper, the problem's author deeply understood both the challenge, and inherent beauty, presented by simplicity.

# Chapter 1

## Introduction

As a light breeze rustles the leaves, listening to the responding clank of wooden chimes recalls the words of Charles Moore: "…breezes loosely captured can connect us with the very edge of the infinite." [1] In our current culture of bright lights dispelling all darkness, it can be a challenge to imagine Japanese culture of more than two centuries ago – one that embraced shadows. Isolated from the rest of the world, Japan found peace; and in that peace, all aspects of its culture flourished.

Transporting ourselves to Japan of 1750, we find an intellectual environment focused on the beauty of simplicity. Walking down the dirt road, we approach the shrine at Gion. We are dazzled by the sunlight dancing off its bright red paint, amazed at the shrine's fascinating angles and symmetry. Below the shrine's eaves, the wooden tablets and paper lanterns flutter in the light wind; handcrafted objects made as offerings to Shintō *kami*, the deities of Japan's native religion. Walking under the long extended roof, entering the shrine produces a feeling of anticipation and depth. There is mystery in every corner - the suspense found in darkness as it enhances the light. "…we find beauty not in the thing itself but in the patterns of shadows, the light and the darkness, that one thing against another creates." [2] Lit only by candles, and light filtering in through periodic paper panels, the shrine is filled with imaginary whispers created by dimness. "The quality that we call beauty,.. must always grow from … our ancestors forced to live in dark rooms, [who] presently came to discover beauty in shadows; ultimately to guide shadows towards beauty's ends." [3]  The flicker of the candles is captured and enhanced by the highly polished lacquer and hints of gold in the surroundings.

---

[1] J. Tanizaki, *In Praise of Shadows,* Leete's Island Books, CT 1977, Forward
[2] J. Tanizaki, *In Praise of Shadows,* Leete's Island Books, CT 1977, p 30

"The small serves to illustrate beauty as well as, perhaps better than, the great; and the poet takes a cultivated delight in simplicity." [4]  Delight is found in the simple aspects of life; in finding the simple forms that exist in complexity.

This is the world in which the Gion Shrine problem was created; a world of beauty found in simplicity; and pleasure found in the presence of shadows. It was a time of peace and prosperity.

---

[3] J. Tanizaki, *In Praise of Shadows,* Leete's Island Books, CT  1977, p 18
[4] W. S. Morton, J. K. Olenik, *Japan: Its History and Culture*, McGraw-Hill Professional; New York, 2005,  p 130

# Chapter 2

## Historical Background

As with any aspect of the human condition, nothing happens in isolation. Discussing the Gion Shrine problem, and offering an alternative solution, cannot be done in the absence of an overall historical perspective. Japanese mathematics of the time, now called *wasan,* and the creation of *sangaku*, happened as an extension and part of a unique period in Japanese history.

### 2.1     Edo Period Japan

Constant turmoil is an apt description of sixteenth century Japan. Referred to as an Age of Warring States, *sengoku*,[5] Japan was ruled by a multitude of lords, *daimyōs*. Fierce competition occurred for resources and attention of the Shōgun (leader of the Japanese military and government), resulting in constant internal warfare. Little sense of a central government or of a unified nation existed.

With Tokugawa Ieyasu winning the Battle of Sekigahara in 1600, change began. Appointment of Tokugawa as Shōgun in 1603, created the reign of the Tokugawa Shōgunate, allowing Ieyasu and his descendants to propagate a peace that lasted for nearly 300 years. In this peace, all aspects of Japanese culture flourished.

---

[5] "Sengoku period", Wikipedia, https://en.wikipedia.org/wiki/Sengoku_period

Beginning in the mid-sixteenth century, trade with Europe brought not only financial gain but also Christian missionaries. Upon becoming Shōgun, Tokugawa Ieyasu implemented a ban against Christianity, forcing the religion to go underground. By continuing this religious prohibition, Ieyasu's descendants[6] encouraged Confucianism to rise, joining Shintō and Buddhism as Japan's religions. Confucianism's ascent changed the Japanese culture in significant ways. Confucianists subscribe to high levels of education and self-improvement, resulting in increases of educational opportunities throughout Japan. Confucianism encouraged reigning in turbulence and replacing it with calm; it embraced straightforwardness and frankness. Japan's literature and art responded with simplicity fostering directness.

Trade with foreign countries initially continued, as did travel overseas. Becoming shōgun in 1623, Iemitsu, Ieyasu's grandson, initiated the period of isolation, *sakoku*. Iemitsu reduced trade to be with the Dutch only. Travel overseas was forbidden as was contact with other Asian countries. Japan became an isolated country.

Under Iemitsu, *daimyōs* were reduced in number through consolidation instead of conquest.[7] Building an administrative structure in the government added administrative positions. A police state was created enforcing laws uniformly, resulting in stability across the nation. Implemented taxation laws focused on agricultural success, increasing the wealth of the lower social classes.

---

[6] H. Fukagawa, T. Rothman, *Sacred Mathematics: Japanese Temple Geometry;* Princeton University Press, Princeton, N.J. 2008 p 4

[7] *Companion to Japanese History,* Blackwell Publishing, 2007; P. C. Brown, *Chapter 4: Unification, Consolidation and Tokugawa Rule;* pp 70

Upon becoming Shōgun, Ieyasu moved Japan's capital from Kyoto to Edo. Enforced capital visitation was imposed requiring *daimyōs* to spend every other year in the capital city. In years when the *daimyō* was not a city resident, he was required to leave family members still in the capital. Forced visitations impacted Japan from multiple perspectives. Travel increased, as did urban development. Thoroughfares required inns, taverns, souvenir shops, etc. to support the traveling *daimyōs*. Increased travel initiated more road maintenance and bridge repair. Required to support two houses, *daimyōs* were no longer able to support their samurai troops. The overall growth in prosperity was felt throughout Japanese society.

The Japanese began to view the Shōgun as being "wrapped in virtue." [8]

## 2.2     Intellectual Aspects of Edo Period Japan

*Sakoku* was a peaceful time where the samurai put away their swords. In 1615, Ieyasun issued "Ordinances for Military Families".[9] Samurai were beseeched to give equal time to both military practices and refined pursuits. Encouraged to read the Japanese classics, learn medicine, gain musical skills, increase mathematics and calculation knowledge, the samurai class changed. With salaries based on the fluctuating price of rice, samurai were forced to find additional sources of income.

Finding administrative positions, samurai brought with them their attitudes of disciple and devotion, altering the governmental climate and attitude.

---

[8] *Companion to Japanese History,* Blackwell Publishing, 2007; P. Nosco *Chapter 6: Intellectual Change in Tokugawa's Japan;* pp 110
[9] W. S. Morton, J. K. Olenik, *Japan: Its History and Culture*; McGraw-Hill Professional; New York, 2005, pp 127

Samurai founded military academies, became artists, literary figures, and mathematicians. Matsuō Bashō (1644-1694), one of the greatest haiku poets who has ever lived, was from samurai descent.

In the absence of universities, many samurai opened schools that taught reading, writing and arithmetic. With certification, samurai opened mathematical academies. Initially restricted to other samurai only, the mathematics academies were later opened to citizens of all class levels. As urbanization and wealth grew, merchants and artisans earned new status. Their burgeoning businesses required these tradesmen to increase their levels of education overall, including mathematics. As the Tokugawa Period continued, mathematical education cut across social levels.

Education actually became a commodity of value that could be compared to a beautiful brocade.[10]

## 2.3    Edo Period Literature and Arts

All of the arts flourished during the Edo Period. More leisure time and higher levels of education, increased the demand for literature of all types. This was particularly true with a form of poetry called haiku. Written by Japanese in a single column, haiku uses "cut words" to create pauses when the poem is read. Although there was some variation, the structure was to have a syllable structure of 5-7-5. This means that the first phrase should have 5 syllables when read, the second phrase, 7 syllables, and the last phrase is to contain 5 syllables. As English is a more condensed language,[11] translated haiku's exact syllable structure is often lost. A three line structure is used in English haiku. The first line usually gives an image. The second line often seems to not completely relate to the first line. The third line then joins the first two lines into a single, unified expression. Left purposefully open-ended, the haiku artist depends on the reader's imagination to complete the poem.

---

[10] *Companion to Japanese History,* Blackwell Publishing, 2007; P. Nosco *Chapter 6: Intellectual Change in Tokugawa's Japan;* pp 111
[11] S. Addiss, *The Art of Haiku;* Shambala Publications, Inc., MA 2012 pp 1-2

Creating beautiful haiku requires a precise mind that can understand the essential, simple aspects of life. "…extreme brevity demanded a bright wit…" [12]

Economic reasons often motivated poets to develop painting skills. As wealth increased along all social levels, the arts gained more patrons. Forced capital visitation imposed on *daimyōs*, required them to maintain two houses, each house requiring decorations. Demand for business decorations also increased rapidly.

Combining haiku and painting was a concept that began with the Chinese. The historical mindedness of the Edo period, motivated some haiku poets to embrace Chinese verse and painting. One of these was the master haiku poet Matsuō Bashō (1644-1694). Bashō fused his haiku and paintings on the same scroll, creating stunning images. The combination of painted images accompanying haiku is called *haiga*. Figure 1 shows a scroll created by another master haiku poet and *haiga* artist, Yosa Buson (1716-1784).

---

[12] W. S. Morton, J. K. Olenik, *Japan: Its History and Culture*; McGraw-Hill Professional; New York, 2005, pp 129

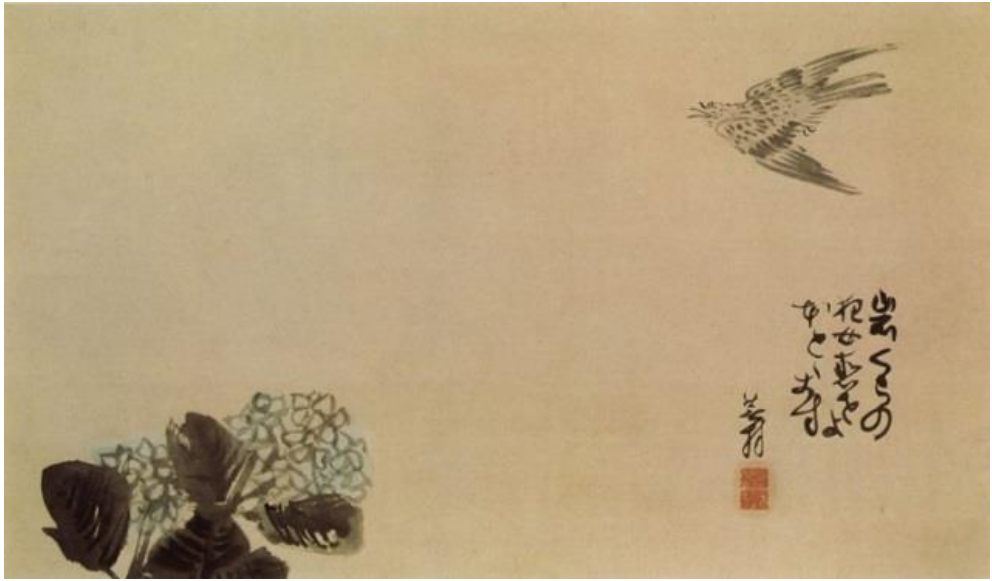

**Figure 1: "A little cuckoo across a hydrangea", by Yosa Buson [24]**

### 2.4 *Wasan* – Japanese Mathematics

*San*, literally meaning arithmetic or counting, was the term used during the Edo Period for the practice of mathematics. Once *sakoku* ended, the term *wasan* was created to indicate the mathematics performed during the Edo Period.[13] *Wasan* practitioners are referred to as *wasanka*.[14] For most of the Tokugawa Shōgunate, the center of mathematics was in Kyoto. It was not until the nineteenth century, after a peasant revolt in 1805, that the center of mathematics moved to Edo, the capital at that time.[15]

---

[13] J. Shigeru, (H. Selin, editor) *Mathematics Across Cultures: The history of Non-Western Mathematics (The Dawn of Wasan;* Kluwer Academic Publisher; Dordrecht, 2000 p 423
[14] J.M. Unger, *Sangaku Proofs: A Japanese Mathematician at Work;* Central East Asia Program, NY 2015, p *x*
[15] J. Shigeru, (H. Selin, editor) *Mathematics Across Cultures: The history of Non-Western Mathematics (The Dawn of Wasan;* Kluwer Academic Publisher; Dordrecht, 2000 p 429

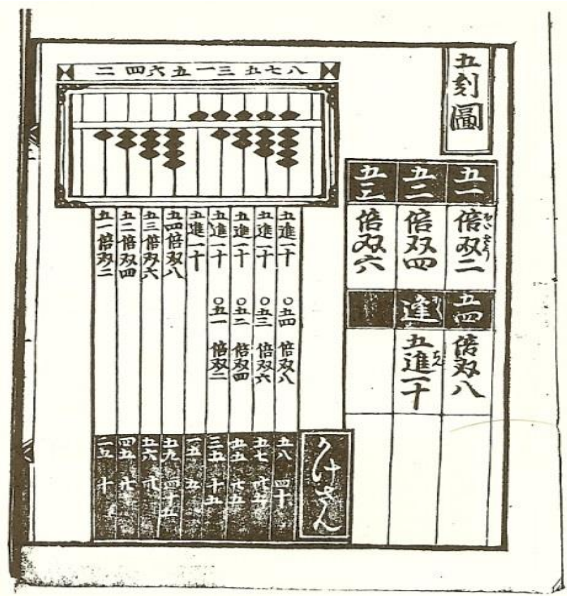

**Figure 2: Picture of *soroban* in upper left corner from *Jinko-ki* [25]**

Most pre-Edo Period mathematics in Japan was based on Chinese texts. Greatly influential to both Chinese and Japanese mathematicians, Cheng Da-Wei's *Systemic Treatise on Arithmetic,* published in 1593, entered Japan just prior to *sakoku*. Cheng's *Treatise* influenced Yoshida Mitsuyoshi's widely used mathematics text *Jinko-ki*, (printed in 1643), a page of which is shown in Figure 2. Yoshida's book covers the use of the *soroban*, the Japanese abacus, a tool that grew rapidly in popularity during *sakoku.* In 1676, Takebe Katahiro published a Japanese rendition of Cheng's *Treatise*.[16]

---

[16] H. Fukagawa, T. Rothman, *Sacred Mathematics: Japanese Temple Geometry;* Princeton University Press, Princeton, N.J. 2008 p 30-31

Edo Period *wasanka* knew Euclidean geometry. Both the quadratic equation and Pythagorean theorem were often utilized. Japanese knew squares, multiplicity of roots,[17] along with trigonometry and logarithms.[18] Pi was known. In fact, several Japanese mathematicians during the Edo Period, calculated new values of pi with great levels of correctness and numerous digits. Fractions were freely used although decimals were also known.[19] Two linear equations could be solved. *Wasanka* knew how to determine area and volume. As discussed later, many of these techniques were used when solving *sangaku*.

Geometric relationships were freely studied in Japan, but most studies were algebra based.[20] Prior to the Edo Period, mathematics was viewed as a branch of science. S*akoku* cultivated a shift toward viewing mathematics as an art. The impact created by visualizing mathematics as an art was one of speciality. Moving away from generality, mathematics became focused on lengthy calculations. During the eighteenth century, *wasanka* discussed returning to more generalized approaches to mathematics as a whole. "The writers of the 18th and 19th centuries had here and there tried some problems by the geometrical way, however imperfect they might have been." [21] Reviewing Ajima Naonobu's work, experts state that a considerable amount of generalization was accomplished by Ajima's efforts.[22]

---

Although mathematical knowledge in Japan expanded during the Edo Period, the country did not have access to the simultaneous mathematical developments in Europe. These accomplishments included calculus. Modern *sangaku* solutions that involve calculus could not have been solutions intended by *wasanka* creating *sangaku.*

## 2.5 *Sangaku*

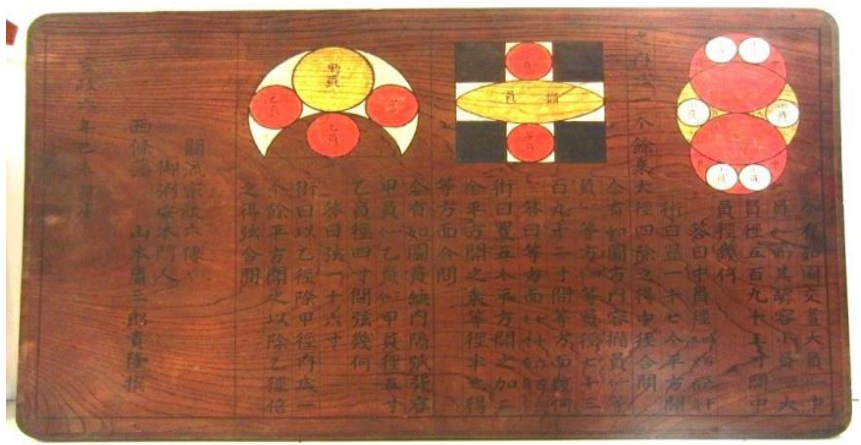

**Figure 3:** *Sangaku* **of Konnoh Hachimanga Shrine, Shibuya, Japan (1859) [26]**

All of *wasan* was viewed as an art during *sakoku.* This is epitomized by *sangaku*, hung in shrines and temples as offerings to the various gods, or "votive plaques".[23] S*an* is defined as arithmetic or calculations. The term *gaku* means plaque,[24] with the literal interpretation of *sangaku* being arithmetic or calculation plaque.

Over the past couple of decades, debate has occurred as to the intention of *sangak*u. During the 1990s, it appears that the primary modern viewpoint concerning *sangaku*'s purpose was one of a religious offering. This viewpoint shifted as the new millennium approached.[25] The 21st Century shift was towards seeing *sangaku* more as a "bragging" or competition tablet, where the author demonstrated his superior mathematical talent.

---

[23] J.M. Unger, *Sangaku Proofs: A Japanese Mathematician at Work;* Central East Asia Program, NY 2015, p *xiv*
[24] R. J. Hosking, *Solving Sangaku: A Traditional Solution to a Nineteenth Century Japanese Temple Problem,* https://arxiv.org/ftp/arxiv/papers/1702/1702.01350.pdf

I propose that *sangaku* were originally votive tablets to the *kami*. As *sakoku* progressed, *sangaku* became used as a competitive tool. This is based on the changing nature of the Japanese culture during *sakoku*. At the beginning of *sakoku*, Japan was just coming out of a period of constant warring, *sengoku*, and finding peace after a century of constant internal chaos. In the early 1600s, education and wealth was primarily restricted to the upper classes. As *sakoku* progressed, education increased across classes, even becoming a commodity about which one could brag. Wealth also spread across the classes. *Sakoku* of the 1700s and 1800s had a very different "personality' than just a century earlier. *Sangaku* can be viewed as representing many aspects of the entire Edo Period, including social and cultural significance.[26]

It is estimated that about 900 *sangaku* exist today. Beginning with the initial documentation, indications show that several thousand *sangaku* were created.[27] The primary reason so few *sangaku* have survived is that the tablets are wooden and wood deteriorates over time. Another reason is the volcanic nature of Japan's island. Japan not only has volcanic eruptions (the most recent one being in 2014), the island also has earthquakes and even periodic tsunami.

---

[25] R. Hosking, *Sangaku: A Mathematical, Artistic, Religious, and Diagrammatic Examination,* Univ. of Canterbury, UK, PhD. Thesis *https://ir.canterbury.ac.nz/handle/10092/12912*  p 3

[26] R. Hosking, *Sangaku: A Mathematical, Artistic, Religious, and Diagrammatic Examination,* Univ. of Canterbury, UK, PhD. Thesis *https://ir.canterbury.ac.nz/handle/10092/12912*  p 6, 257

[27] R. Hosking, *Sangaku: A Mathematical, Artistic, Religious, and Diagrammatic Examination,* Univ. of Canterbury, UK, PhD. Thesis *https://ir.canterbury.ac.nz/handle/10092/12912*  p 24

Although all *sangaku* were geometric in nature, the variation of problem content written on them is significant. Some *sangaku* posed agricultural problems having to do with harvests and field sizes while others were abstract geometric shapes, requiring challenging solutions. The oldest known *sangaku* dates to 1668, with the first published collection of *sangaku* authored by Fujita Kagen (1765-1821). Published in 1789, Kagen's collection was titled *Mathematical Problems Suspended Before The Temple.*[28] Although most *sangaku* were created by the members of the samurai class, research on *sangaku* reflects that people from all social levels authored the tablets.[29]

*Sangaku* painted images did not always exactly match the descriptions given for the problem being presented.[30] This is more than likely due to the fact that the *sangaku* problem author was not necessarily the same person who painted the *sangaku* image. If they did not feel capable of painting a *sangaku* design, *wasanka* would hire an artist to do so.

There are well known instances where the creator of both the *sangaku* mathematical problem and the *sangaku* image were the same individual. One *sangaku* pictured in Rosalie Hosking's paper, *Sangaku: A Mathematical, Artistic, Religious, and Diagrammatic Examination*, the Yoshifuji Mishima shrine tablet, is known to have been authored by a painter. "It is known that this particular mathematician was skilled in many

---

[28] C.A. Pickover, *The Math Book;* Sterling, NY 2001 p 198
[29] T. Rothman, H. Fukagawa, Japanese Temple Geometry, *Scientific American,* Vol. 278, No. 5 (May, 1998, pp. 84-91)
[30] R. Hosking, *Sangaku: A Mathematical, Artistic, Religious, and Diagrammatic Examination,* Univ. of Canterbury, UK, PhD. Thesis *https://ir.canterbury.ac.nz/handle/10092/12912* p 266

areas, and alongside the *sangaku* hangs a painting of Mount Fuji painted by the same author." [31]

Created by members of all social classes, *sangaku* were also intended to be viewed by everyone. In his travel diary, mathematician Yamaguchi Kanzan (1781-1850), mentions that two other mathematicians tell him "…we decided to hang a *sangaku* in this shrine. We hope that the visitors will look at this tablet and ask for any opinions about the problem." [32]

Figure 4 displays reproductions of six *sangaku* geometric images. Although some *sangaku* geometric designs were complex, the designs displayed in Figure 4 were purposefully chosen to be less complex images.

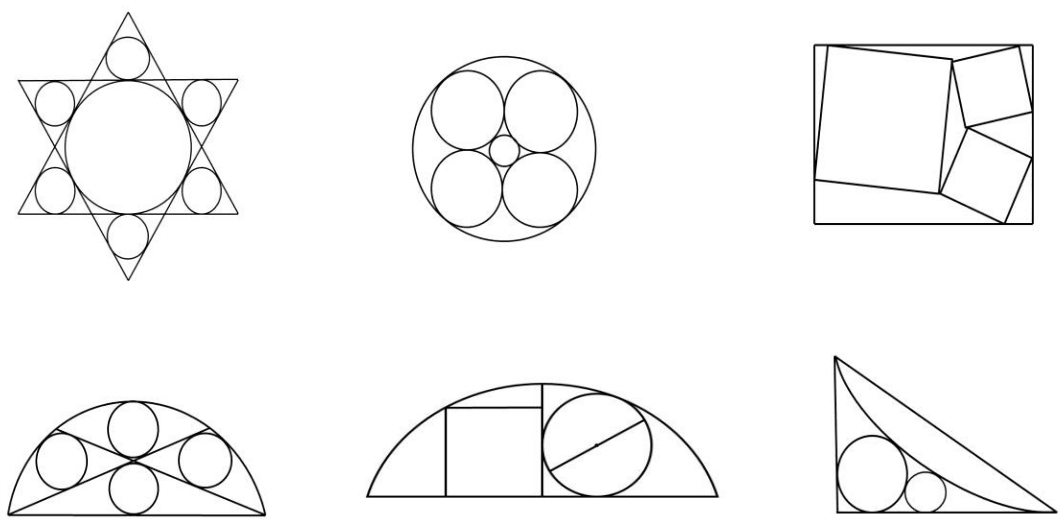

**Figure 4: Six *sangaku* examples.**

---

[31] R. Hosking, *Sangaku: A Mathematical, Artistic, Religious, and Diagrammatic Examination,* Univ. of Canterbury, UK, PhD. Thesis *https://ir.canterbury.ac.nz/handle/10092/12912* p 251

[32] R. Hosking, *Sangaku: A Mathematical, Artistic, Religious, and Diagrammatic Examination,* Univ. of Canterbury, UK, PhD. Thesis *https://ir.canterbury.ac.nz/handle/10092/12912* p 253-254

The Gion Shrine geometric design is displayed in the middle of Figure 4's second row. Although the Gion Shrine design is simple, and in context similar to that of the others, it has an asymmetrical balance and grace that sets it apart.

Labels in the *sangaku* figure and relevant unknowns were identified using characters from the Chinese calendar.[33] Today's mathematics uses characters from both the English and ancient Greek alphabets in a manner similar to that of calendar characters on *sangaku*.

A typical *sangaku* had three parts.[34] These three aspects will be used in section 3.2, "The Gion Shrine *Sangaku*", to analyze four modern presentations of the Gion Shrine *sangaku*.

1. The Problem: The problem statement consisted of an image often with accompanying labels and descriptive text. The written portion sometimes provided numerical values for measures on the *sangaku* image, but this was not always the case.

2. The Answer: Often a numerical answer was given. Full worked out solutions were rare; while some *sangaku* simply referred to the section with the formulae, not offering an answer.

3. Formulae: This section gave the information, technique or formulae, required to solve the problem.

---

[33] R. J. Hosking, *Solving Sangaku: A Traditional Solution to a Nineteenth Century Japanese Temple Problem,* https://arxiv.org/ftp/arxiv/papers/1702/1702.01350.pdf

[34] R. J. Hosking, *Solving Sangaku: A Traditional Solution to a Nineteenth Century Japanese Temple Problem,* https://arxiv.org/ftp/arxiv/papers/1702/1702.01350.pdf

"To learn traditional Japanese mathematics is to learn another way of thinking." [35] Going through the *sangaku* solutions offered by various authors shows this statement to be true. J. M. Unger [36] presents twenty-six *sangaku* solutions in *Sangaku Proofs: A Japanese Mathematician at Work*. Rosalie Hosking in her paper *Sangaku: A Mathematical, Artistic, Religious, and Diagrammatic Examination*, offers numerous *sangaku* with very detailed translations and complete solutions. Hosking's paper also includes photos of existing *sangaku*, adding a delightful dimension to her work. H. Fukagawa and T. Rothman, [37] in *Sacred Mathematics: Japanese Temple Geometry*, have solutions located separately from the problems, challenging readers to try working the *sangaku* problems themselves.

## 2.6    *Sangaku* Translation

Reflecting on the historical mindedness of *sakoku*, *sangaku* were written in a cryptic language, *kanbun*, closely related to ancient literary Chinese. Other than diagrammatical marks indicating that it is Japanese, *kanbun* is essentially Chinese characters following Chinese grammar. *Kanbun* can be compared to the use of Latin in some older European mathematical texts, such as those written by Newton and Euler. *C*onsidered by Japanese society as an indication of higher levels of education, Samurai were taught *kanbun* at an early age. Rosalie Hosking mentions the following in reference to the language used on *sangaku*: "The language of *sangaku* is quite conceptual, with often just the core nouns and verbs present…. Because of its conceptual nature, the language of *sangaku* can come across as vague."[38] A different more explicit language was used in texts written for merchants, famers and artisans.[39]

---

[35] H. Fukagawa, T. Rothman, *Sacred Mathematics: Japanese Temple Geometry;* Princeton University Press, Princeton, N.J. 2008 p 2

[36] J.M. Unger, *Sangaku Proofs: A Japanese Mathematician at Work;* Central East Asia Program, NY 2015

[37] H. Fukagawa, T. Rothman, *Sacred Mathematics: Japanese Temple Geometry;* Princeton University Press, Princeton, N.J. 2008

In his book, Unger also discusses *wasanka*'s obsessiveness over writing *sangaku* problems in as few characters as possible.[40] This fits well with *sakoku* where succinctness was valued not only in writing mathematics but also in literary works such as haiku.

Table 1 shows a comparison taken directly from Rosalie Hosking's *Sangaku: A Mathematical, Artistic, Religious, and Diagrammatic Examination*. The Japanese characters accompanying the text in Hosking's paper are omitted here, creating yet another "translation error"! Table 1 not only displays differences between Hosking's literal translation versus a modern translation, it allows the reader of this paper to better understand the cryptic style of *kanbun* as shown in Hosking's direct translation.

**Table 1: Comparison of Sangaku Translations**

| | *Literal Translation by Hosking* | *Modern Translation by Nagano Wasan Society (2005)* |
|---|---|---|
| Problem: | There is a diagram line five circle on say *ko* diameter 2 *sun otsu* circle diameter 3 *sun tei* circle diameter how much. | The circles *ko*, *otsu*, and *tei* are neighboring and circumscribe the line. Also, four other circles circumscribe circle *hei*. When the circles *ko* and *otsu* are 2 sun and 3 sun respectively, find the diameter of circle *tei*. |
| Answer: | *tei* circle 18 *sun* | Circle *tei* diameter 18 sun |
| Formulae: | Technique by means of *ko* circle diameter 2 and *otsu* circle diameter subtract *otsu* circle diameter divide square *ko* circle multiply obtain. | Using circles *ko*, *otsu*, *hei* and *tei* respectively as $d_1, d_2, d_3,$ , and $d_4$. $d_4 = \frac{(d_1 \cdot d_2)^2}{(2d_1 - d_2)^2}$ |

When comparing the two translations in her paper, Hosking says the following:[41]

"The creators of *sangaku*, though having knowledge of symbolic manipulation through *tenzan jutsu…* specifically chose to present solutions verbally in the technique section. By replacing this text with modern algebraic formulas, the text no longer accurately represents the original problem and intention of the author."  It should be noted that both translations allowed Hosking and the Nagano Wasan Society to reach the same answer.

---

[41] R. Hosking, *Sangaku: A Mathematical, Artistic, Religious, and Diagrammatic Examination,* Univ. of Canterbury, UK, PhD. Thesis *https://ir.canterbury.ac.nz/handle/10092/12912* p 12-13

# Chapter 3

# The Gion Shrine Problem

One *sangaku*, called the Gion Shrine problem, has evaded solution for at least years.[42] The original Gion Shrine *sangaku* is long gone due to the passage of time. The original *sangaku* went undocumented, with its author unidentified.

### 3.1    The Gion Shrine *Sangaku* History

In 1749, Tsuda Nobuhisa posted the first solution *sangaku* to the Gion Shrine problem.[43] The solution Tsuda offered had a variable degree of 1024.[44] This solution was followed by another solution given by the mathematician Nakata that reduced the degree to 46. A third solution was given by Ajima Naonobu in 1774. Becoming well known in his lifetime partly for his solution to the Gion Shrine *sangaku*, Ajima published *Kyoto Gion Gaku Toujyutsu*, giving an answer reflecting a tenth degree polynomial.[45] Ajima received such high praise for his accomplishment that he also received a high level position in the government.[46] The technique that Ajima utilized was based on simple geometric applications.[47] There is mention of a fourth *sangaku* rewriting done by Saito Mitsukuni in 1815.[48]

---

[42] H. Fukagawa, T. Rothman, *Sacred Mathematics: Japanese Temple Geometry;* Princeton University Press, Princeton, N.J. 2008 p 250

[43] H. Fukagawa, T. Rothman, *Sacred Mathematics: Japanese Temple Geometry;* Princeton Univ. Press, Princeton, N.J. 2008 p 250

[44] Y. Mikami, *The Development of Mathematics In China and Japan;* Chelsea Publishing Company, UK, 1974 pp 168

[45] J. A. de Reyna, D. Clark, N. D. Elkies, *A Modern Solution to the Gion Shrine Problem*, https://folios.rmc.edu/davidclark/wp-content/uploads/sites/56/2016/12/gion.pdf

[46] J.J. O'Connor, E.F. Robertson, *Chokuyen Naonobu Ajima;* https://www-history.mcs.st-andrews.ac.uk/Biographies/Ajima.html

Modern solutions have involved polynomials requiring mathematics not available to Japan at the time the Gion Shrine problem created. For this reason, the Gion Shrine problem is still considered unsolved.[49]

A modern solution for the Gion Shrine *sangaku* is presented by de Reyna, Clark and Elkies in their paper, *A Modern Solution to the Gion Shrine Problem*. In addition to giving proof for their Gion Shrine *sangaku* solution, de Reyna, Clark, and Elkies also prove that rational solutions to this problem are not possible.[50] I agree with this conclusion and propose that the Gion Shrine problem author intended his *sangaku* to not to be solvable by standard means.

### 3.2    The Gion Shrine *Sangaku*

With its beautiful balance and elegance, the Gion Shrine *sangaku* has received significant attention. For more than  years, people from diverse walks of life have offered solutions.

---

[47] J. A. de Reyna, D. Clark, N. D. Elkies, *A Modern Solution to the Gion Shrine Problem*, https://folios.rmc.edu/davidclark/wp-content/uploads/sites/56/2016/12/gion.pdf
[48] H. Fukagawa, T. Rothman, *Sacred Mathematics: Japanese Temple Geometry;* Princeton Univ. Press, Princeton, N.J. 2008 p 250
[49] H. Fukagawa, T. Rothman, *Sacred Mathematics: Japanese Temple Geometry;* Princeton University Press, Princeton, N.J. 2008, p 250
[50] J. A. de Reyna, D. Clark, N. D. Elkies, *A Modern Solution to the Gion Shrine Problem*, https://folios.rmc.edu/davidclark/wp-content/uploads/sites/56/2016/12/gion.pdf

### 3.2.1 Modern Interpretations: The Problem – The Design

Modern documented images of the Gion Shrine *sangaku* vary as does the statement of the *sangaku*'s problem description and statement of the formulae. Figure 5 reflects four reproductions of modern images that have been published for the Gion Shrine *sangaku*.

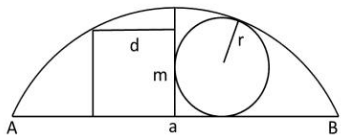

H. Fukagawa, T. Rothman, *Sacred Mathematics: Japanese Temple Geometry* [2]

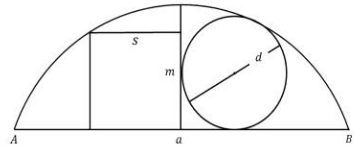

J.A. de Reyna, D. Clark, *A Modern Solution To The Gion Shrine Problem* [5]

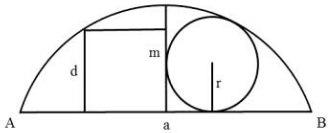

Cut-the-knot; Gion Shrine Problem; http://www.cut-the-knot.org/pythagoras/GionShrineProblem.shtml [27]

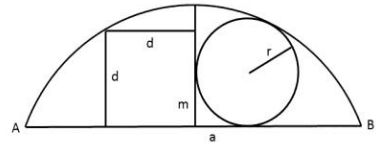

J.J. O'Connor, E.F. Robertson, *Chokuyen Naonobu Ajima* [6]

**Figure 5: Four current images of Gion Shrine *sangaku* design.**

To create the reproductions in Figure 5 as closely as possible to each modern rendition of the Gion Shrine *sangaku* design, a strict procedure was followed. Copies of each modern image were obtained from either the printed documents or the web. These copies were pasted into PowerPoint. Overlapping the pasted copies, the reproductions were created matching the copies as accurately as possible. 200% zoom was used to increase accuracy of each reproduction to its matching modern day image. Once the reproduction was completed, the modern copy was deleted from PowerPoint. Slight differences might be found in fonts used for labels and label placement on the reproduction.

Notice the differences between the images displayed in Figure 5. There are differences in curvature, labels used and label placement. All of these images have a capital "A" and a capital "B" label, although these are only mentioned in the problem's description when defining the length of chord $a$.

Notice that three of the images in Figure 5 have radii and one image has a diameter with a clear center point. Also note that the placement of the radii varies significantly, making a difference in visual beauty and balance. Three different angle measures are represented by the radii and diameter of the images shown in Figure 5.

### 3.2.2 Modern Interpretations: The Problem Description

In the modern presentations of the Gion Shrine *sangaku*, there are also differences in the statements of problem description. Utilizing the same sources as those in Figure 5, following are the authors' problem descriptions.

1. H. Fukagawa and T. Rothman, *Sacred Mathematics: Japanese Temple Geometry*

    "We have a segment of a circle. The line segment $m$ bisects the arc and

chord $AB$. As shown, we draw a square with side $d$ and an inscribed circle of the radius $r$. Let length $AB = a$." [51]

2. J.A. de Reyna and D. Clark, *A Modern Solution To The Gion Shrine Problem*

   "We have a segment of a circle. The line segment $m$ bisects the arc and chord $AB$. As shown, we draw a square with side $s$ and an inscribed circle of diameter $d$. Let the length $AB = a$." [52]

3. Cut-the-Knot web page; "Gion Shrine Problem"

   "In a circular segment with base $AB$ of length $a$ and altitude $m$, there are a circle of radius $r$ inscribed in one half of the segment and a square of side $d$ inscribed in the other half, as shown." [53]

4. J.J. O'Connor, E.F. Robertson, *Chokuyen Naonobu Ajima*

   "In this figure we have a segment of a circle on the chord $AB$ of length $a$. From the mid-point of $AB$ we draw a line perpendicular to $AB$ to meet the circle. It has length $m$. To the left of this line we draw a square of side $d$, as shown, and to the right we draw a circle of radius $r$, as shown." [54]

In the four expressions of the Gion Shrine problem, the descriptions given by Fukagawa/Rothman and de Reyna/Clark are almost identical. De Reyna and Clark use different variables that coincide with the fact that their image displays a diameter instead of a radius. Both of the first two modern descriptions include "inscribed", a term not utilized by *wasanka*. Cut-the-Knot web site utilizes two terms, "inscribed" and "altitude", that are not ones Edo Period Japanese would have used. The description offered by O'Connor and Robertson reflects the design's construction. As evidenced by the direct translation in column one of Table 1, this approach is typical of the ones used by *wasanka* in describing *sangaku* images.

---

[51] H. Fukagawa, T. Rothman, *Sacred Mathematics: Japanese Temple Geometry;* Princeton University Press, Princeton, N.J. 2008 p 250
[52] J.A. de Reyna, D. Clark, *A Modern Solution To The Gion Shrine Problem;* (2013) https://arxiv.org/abs/1306.5339
[53] Cut-the-knot; Gion Shrine Problem; http://www.cut-the-knot.org/pythagoras/GionShrineProblem.shtml
[54] J.J. O'Connor, E.F. Robertson, *Chokuyen Naonobu Ajima;* https://www-history.mcs.st-andrews.ac.uk/Biographies/Ajima.html

In their book, Fukagawa and Rothman state that they used the rewritten version of the Gion Shrine problem posted by Saito Mitsukuni in 1815. O'Connor and Robertson state in their article that they utilized Ajima's rewriting in his *Kyoto Gion Dai Toujyutsu* published in 1774. The variation between the words used by O'Connor and Robertson versus those chosen by Fukagawa and Rothman might not appear to be significant, but I contend that the differences are important. Also consider the fact that in the few years between 1774 and 1815 (only 41 years), the wording and the image of the Gion Shrine problem appears to have changed as presented by modern interpretations.

### 3.2.3  Modern Interpretations: The Answer

The first solution *sangaku* to the Gion Shrine problem was posted  years ago. Including Tsuda's initial solution *sangaku*, no one offering a solution has mentioned that the original *sangaku*, hung prior to 1749, included a numerical answer. This includes the four modern sources being compared in this paper. It can be safely concluded that no answer was given on the original tablet.

### 3.2.4 Modern Interpretations: The Formulae

In the four modern presentations of the Gion Shrine *sangaku*, there are also differences in the statement of the formulae. The modern sources discussed previously are compared following.

5. H. Fukagawa and T. Rothman, *Sacred Mathematics: Japanese Temple Geometry*

    "Then, if $p = a + m + d + r$ and $q = \frac{m}{a} + \frac{r}{m} + \frac{d}{r}$, find $a, m, d,$ and $r$ in terms of $p$ and $q$." [55]

6. J.A. de Reyna and D. Clark, *A Modern Solution To The Gion Shrine Problem*

    "Then, if $p = a + m + s + d$ and $q = \frac{m}{a} + \frac{d}{m} + \frac{s}{d}$, find $a, m, s,$ and $d$ in terms of $p$ and $q$." [56]

7. Cut-the-Knot web page; Gion Shrine Problem

    "Form $p = a + m + d + r$ and $q = \frac{m}{a} + \frac{r}{m} + \frac{d}{r}$. The task is to express $a, m, d,$ and $r$ in terms of $p$ and $q$." [57]

8. J.J. O'Connor, E.F. Robertson, *Chokuyen Naonobu Ajima*

    "Put $p = a + m + d + r$, and $q = \frac{m}{a} + \frac{r}{m} + \frac{d}{r}$. The problem requires that we express $a, m, d,$ and $r$ in terms of $p$ and $q$." [58]

---

[55] H. Fukagawa, T. Rothman, *Sacred Mathematics: Japanese Temple Geometry;* Princeton University Press, Princeton, N.J. 2008 p 250
[56] J.A. de Reyna, D. Clark, *A Modern Solution To The Gion Shrine Problem;* (2013) https://arxiv.org/abs/1306.5339
[57] Cut-the-knot; Gion Shrine Problem; http://www.cut-the-knot.org/pythagoras/GionShrineProblem.shtml

As with the design description, the formulae given by Fukagawa/Rothman and de Reyna/Clark, with the exception of two variables, are the same. There are meaningful differences between Cut-the-Knot's formulae wording from that of Fukagawa and Rothman. Where Fukagawa and Rothman use "Then, if" and "find", Cut-the-Knot use "Form" and "The task is to express". O'Connor/Robertson and Cut-the-Knot have similar wording with both of these interpretations being quite different from the first two sources in the list. One critical difference between sources three and four is the use of "Form" (Cut-the-Knot) versus "Put" (O'Connor and Robertson). Although Cut-the-Knot adds additional words, the web site and O'Connor and Robertson both use the word "express" in their formulae. "Express" fits perfectly with the solution presented in geometry. The wording variation is key. I submit that there are several misinterpretations that have occurred over the years with regards to the Gion Shrine *sangaku*, as will be discussed later in this paper.

Section 3.1, "The Gion Shrine *Sangaku* History", mentions that the first solution to the Gion Shrine *sangaku* was posted by Tsuda Nobuhisa in 1749. The idea that Nobuhisa might have slightly altered the formulae, is a very reasonable assumption to make.

On page 39 of Hosking's paper, she states the following concerning her translation of *sangaku*: "To translate *sangaku*, I have connected *kanbun* characters found on *sangaku* back to related modern Japanese mathematical terminology. For example, the term *ka* which appears on tablets means 'addition' in modern day

---

[58] J.J. O'Connor, E.F. Robertson, *Chokuyen Naonobu Ajima;* https://www-history.mcs.st-andrews.ac.uk/Biographies/Ajima.html

Japanese. It can be assumed that the *sangaku* term *ka* is related, and by treating this character as 'addition' on tablets it can be seen that this is a correct interpretation, for it produces the values on tablets." [59] This statement by Hosking implies that if the assumption is made that *ka* means addition, this interpretation should be checked by making certain that the given values sum to the answer provided. If no numerical values are provided in a *sangaku*'s design description or answer, this verification of interpretation cannot be accomplished.

It can be safely assumed that the original Gion Shrine had no numbers given for any of its variables as none of the translations of this *sangaku* over the past years have offered any. This means that there was no manner in which to check if making the same assumption that Hosking made concerning *ka* (that *ka* indicates addition), was appropriate for the Gion Shrine tablet. However, as seen in the modern interpretations, formulae descriptions for the Gion Shrine *sangaku* include two summations, one for $p$ and one for $q$.

It would also be reasonable to speculate that Tsuda, the first individual to post a Gion Shrine problem solution, might have made the same *ka* assumption, using "+" symbols where none were intended by the original author.

Further in her paper, Hosking says: "While both types [*sangaku* offering numerical values for design parts versus those that do not] essentially have the same types of formula sections (for the formula section is always where one figure is put in terms of others), the exclusion of numerical values in the first type [this would indicate

---

[59] R. Hosking, *Sangaku: A Mathematical, Artistic, Religious, and Diagrammatic Examination,* Univ. of Canterbury, UK, PhD. Thesis *https://ir.canterbury.ac.nz/handle/10092/12912* p 39

*sangaku* such as the Gion Shrine design, that offer no numerical values for image pieces] gives an additional level of abstraction to the problem. For, the answer does not have a numerical value which the reader can check against or work backwards from; all they have to rely on is a deep understanding of the geometry. The language used on *sangaku* thus emphasises understanding the relationships between figures given in the diagram. It is through this understanding that the way to solve the problem could be found." [60]

Without any numerical values given in the Gion Shrine *sangaku*, as Hosking mentions, the diagram becomes essential to the given problem. In fact, it becomes critical to solving the problem. Note the design variations offered by the four modern interpretations of the Gion Shrine *sangaku*.

### 3.2.5 The Gion Shrine *Sangaku*: My Interpretation

My interpretation will follow the structure usually provided by *sangaku* themselves. First, I will focus on my interpretation of the Gion Shrine *sangaku* problem as shown by its design and the accompanying design description. As it has been discussed previously, it is evident that no answer was provided on the original Gion Shrine *sangaku*, so this portion will be omitted. Following my discussion of the Gion Shrine *sangaku* problem, I will offer my version of the formulae for the *sangaku* under study.

---

[60] R. Hosking, *Sangaku: A Mathematical, Artistic, Religious, and Diagrammatic Examination,* Univ. of Canterbury, UK, PhD. Thesis *https://ir.canterbury.ac.nz/handle/10092/12912* p 43

Reviewing the four modern interpretations of the Gion Shrine image as presented earlier in this paper, I conclude that the circle on the right side of the circle segment had a diameter with a well-defined center as shown by de Reyna and Clark.[61] It is known that *sangaku* authors tended to use diameters instead of radii in their designs.[62] The embedded circle image used in this paper's solution shows the angle as given by both de Reyna/Clark and O'Connor/Robertson for the diameter/radius, respectively. I also propose that the slightly broader curvature as shown in some of the images (O'Connor and Robertson notably), more accurately represents the original design.

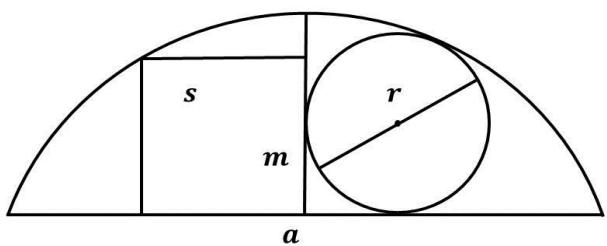

**Figure 6: Gion Shrine *sangaku* geometric design.**

Figure 6 displays the Gion Shrine problem geometric design used in offering the solution in geometry.

The Gion Shrine *sangaku* design was very carefully constructed. I propose that the author carefully constructed this s*angaku* so that only one answer, the answer given in this paper, was possible. The Gion Shrine *sangaku*'s author also did a meticulous

---

[61] J.A. de Reyna, D. Clark, *A Modern Solution To The Gion Shrine Problem;* (2013) https://arxiv.org/abs/1306.5339
[62] J.M. Unger, *Sangaku Proofs: A Japanese Mathematician at Work;* Central East Asia Program, NY 2015, p *xix*

construction so when the correct solution was realized, the solution itself would be as astoundingly beautiful in form as the *sangaku*'s geometric image; both with the magical elegance found in simplicity.

Following is my general description of the Gion Shrine *sangaku* design. It is not meant to be as literal as other modern problem descriptions. In other words, I am not implying that the description I am presenting was on the original Gion Shrine *sangaku*. My intent <u>is</u> to offer a description focused on construction of the geometric image as I believe that correctly matches what would have been on the original *sangaku*. Notice, I am omitting both $A$ and $B$ at the ends of chord $a$ as I believe that they were not included in the original design. With a design focused on simplicity, extraneous labels would not have been included. The only purpose served by $A$ and $B$ is to define the length of chord $a$. This is unnecessary as will be shown.

The Gion Shrine *sangaku* image reflects a circle segment, not a semi-circle. The circle segment contains a perpendicular bisector, $m$, drawn from the bottom line segment (chord), $a$, to the arc above. Vertical line segment $m$ bisects both the arc above and chord $a$ at the bottom of the image. Inside the circle segment in the right half, a circle with a diameter and a well-defined center is drawn. The clear center point in the circle makes it evident that the diameter is $2r$ with $r$ being the radius. The left half contains a square resting on chord $a$. The square touches both the bisector, $m$, and the arc of the circle segment. Inside the circle, is the label $r$. Inside the square is a label $s$. In Figure 6, note that

$$a = 2r + 2m.$$

Although I utilize the diameter as shown by de Reyna and Clark, I conclude that the use of $r$ in the formulae is correct. The formulae closest to fitting with my solution is given following; although I conclude that using a strict mathematical interpretation is not correct.

"Using $p = a + m + s + r$ and $q = \frac{m}{a} + \frac{r}{m} + \frac{s}{r}$,

express $a$, $m$, $s$, and $r$ in terms of $p$ and $q$."

Following, I restate the sentiment of the formulae that I believe was actually written on the *sangaku* containing the design in Figure 6. It is important to keep in mind that *kanbun* is a cryptic language, and one that I do not know. My interpretation here is in modern language with an intent of conveying a general meaning only.

Written on the tablet in *kanbun* characters probably were lines similar to the following: "Put $p$ is $a$ and $m$ and $s$ and $r$" for "$q$ is $m$ related to $a$ and $r$ related to $m$ and $s$ related to $r$"; "Express $a$, $m$, $s$, and $r$ with $p$ and $q$."

It would be very natural for mathematicians, such as ourselves, to read these lines and immediately translate them into mathematical meaning and mathematical symbols, especially if they were written on a *sangaku*.

There is much discussion when it comes to the translation of historical documents. There are two camps with one believing in literal translation and the other believing in making ancient text relevant to today's culture. Either way, as discussed in section 2.6, "*Sangaku* Translation", translation is an issue and I pose it as a possible issue with the Gion Shrine *sangaku*, even when the translation was close to the creation of the original Gion Shrine *sangaku*.

Unger mentions that in translating and working the *sangaku* in his book, he often had to correct "slips of the brush." In 1749, when the first solution to the Gion Shrine sangaku was posted by Tsuda Nobuhisa,[63] I propose that Tsuda felt the same need as Unger, to do corrections. His corrections, which were actually misinterpretations, were then propagated by future mathematicians trying to solve the Gion Shrine *sangaku*.

Individuals with both greater expertise in *kanbun* and with greater access to existing original ancient documents can more accurately judge the truth of these suppositions. However, we are all time's humble servants as the original Gion Shrine *sangaku* no longer exists.

Perhaps the author of the Gion Shrine *sangaku* knew that his words would most likely be misinterpreted by mathematicians. An expert in using words to convey thoughts, the Gion Shrine author was well read in both Chinese verse and *kanbun*.

---

[63] H. Fukagawa, T. Rothman, *Sacred Mathematics: Japanese Temple Geometry;* Princeton Univ. Press, Princeton, N.J. 2008 p 250

# Chapter 4

## The Solution in Geometry

The variation of the images displayed in Figure 5 indicates the difficulty that exists in translating mathematics that belongs in the past. This is especially true when the original *sangaku* no longer exists. Each interpretation is open to each author's own viewpoint. As mathematicians, we bring our entire self to any problem that we address. Although mathematics insists on a high level of precision and accuracy for it to have the greatest meaning, the fact that we each have different viewpoints can be of significant benefit to our field of study. We learned long ago that problems are often solved more completely with the synergy provided by multiple individuals, especially when addressing complex problems. However, variation in interpretations creates uncertainty in knowing truth. A balance between the two situations must be reached. Comparing the Gion Shrine *sangaku* design presented in this paper, Figure 6, and the four modern design interpretations in Figure 5, which image is closest to the original *sangaku* image?

It would be a natural inclination for *wasanka* to assume that the Gion Shrine *sangaku* was like all of the rest. I feel that a *purely* mathematical meaning was not intended and translating the *kanbun* characters on the Gion Shrine *sangaku* into mathematical symbols is making a false assumption that produces a misintended result. The characters painted on the Gion Shrine *sangaku* were to be interpreted just as they were written, not as mathematical statements, requiring mathematical symbols.

My suggested statements for $p$ and $q$ creates a single sentence, "Put $p$ is $a$ and $m$ and $s$ and $r$ for $q$ is $m$ related to $a$ and $r$ related to $m$ and $s$ related to $r$". This is a different manner of interpretation than prior mathematicians. My concept for the solution utilizes $p$ and $q$ to <u>express</u> the various pieces of the *sangaku* design. This emphasizes the last formulae statement, "Express $a$, $m$, $s$, and $r$ with $p$ and $q$.", that I suggest.

The characters used as $p$ and $q$ were meant to represent clues to a mathematical haiku. The line associated with $p$ indicates that the haiku is focused on various parts, <u>not</u> a summation. In this case, it indicates perimeter and circumference of the geometric shapes in the Gion Shrine *sangaku* design.

The formulae line for $q$ shows the relationships that the *sangaku* author meticulously created with his geometric design. The bottom line of the image, $a$, has the relationship where $a = 2r + 2m$. Hence, this is the relationship of $a$ and $m$.

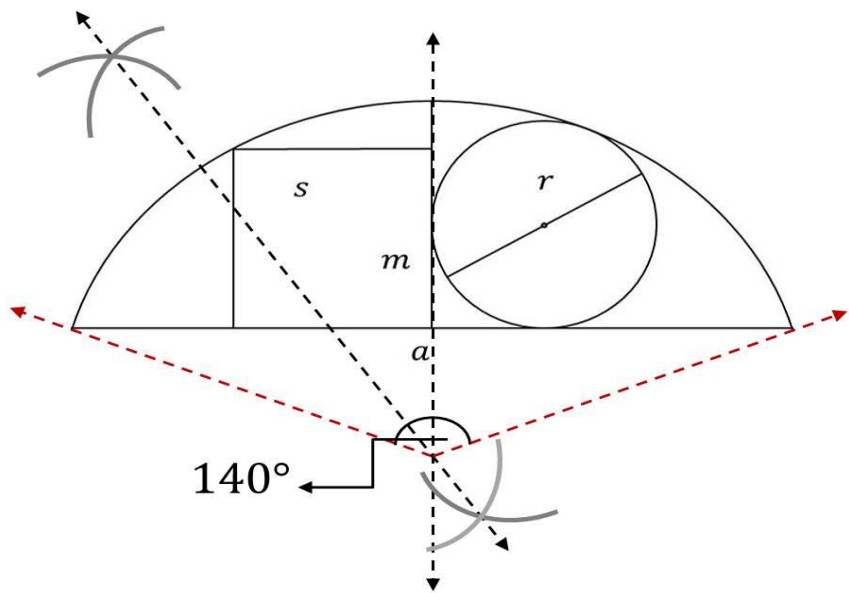

**Figure 7: Gion Shrine geometric design:**
**Locating center point of larger circle and measuring central angle of circle sector.**

The bisecting line segment $m$, when extended downward and intersected with a second line perpendicular to the arc's tangent line, locates the center of the larger circle, part of which is the top of the circle segment. Figure 7 reflects the Gion Shrine *sangaku* image and the calculations required to determine the circle sector's angle.

The bisector $m$ has been extended down. Using two points and the distance between them, the gray arcs represent arcs drawn using a compass as the *wasanka* might have done.[64] A second line perpendicular to the arc's tangent was drawn between the two sets of gray arcs. The center point of the larger circle is where the extended perpendicular bisector intersects the line perpendicular to the arc's tangent. From the center point, lines can then be drawn to the end points of chord $a$ (red dashed lines in Figure 7). The angle between the two red lines is the central angle for the circle sector, part of which is the Gion Shrine *sangaku* design's circle segment.

As *wasanka* would have done, using a protractor, the interior angle of the circle sector in Figure 7 is found to be $140°$. The radius of the large circle, is equal to $2s$. With $a = 2r + 2m$, then the formulae for the perimeter of the circle segment relates $m$ and $r$ as indicated by the second relationship mentioned in statement $q$.

The last relationship given with clue $q$, is "$s$ related to $r$". This relationship indicates the first and the last lines of the mathematical haiku.

Table 2 contains the haiku, in both the *sangaku*'s geometric shapes and the geometry that goes with each shape.

---

[64] R. Hosking, *Sangaku: A Mathematical, Artistic, Religious, and Diagrammatic Examination,* Univ. of Canterbury, UK, PhD. Thesis *https://ir.canterbury.ac.nz/handle/10092/12912*  p 31

**Table 2: Gion Shrine sangaku solution in geometry.**

## *Gion Shrine Perimeter and Circumference Haiku*

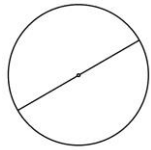          $2r\pi$

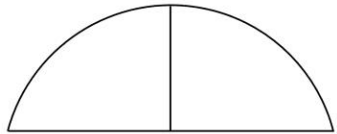          $2r + 2m + \dfrac{7}{9}\pi$

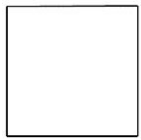          $2s + 2s$

The simple geometric beauty of this solution is as magical as the *sangaku* design.

For reference purposes, the image from Figure 6 is reproduced here.

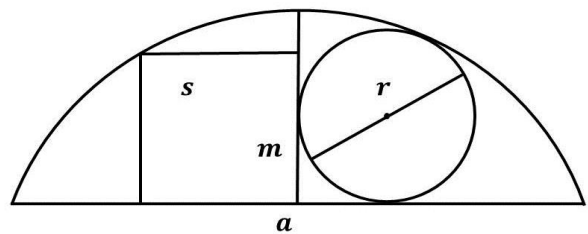

Based on the visual images, note, that in the tradition of haiku, the first line gives an initial image, a circle. The second line seems to relate to the first as it contains a curve; however, as with haiku, the line segments $m$ and $a$ create a question concerning the relationship of the first two lines.  The third line, the square, ties the first two lines together by explaining the straight lines in the second line's image. The original *sangaku* image had the first and last haiku lines contained within the open-ended second line. The outside image wraps around the two inside images and unites them visually. Based on the visual shapes, this completes the *Gion Shrine Perimeter and Circumference Haiku*.

Examining the geometric expressions given, the haiku syllable pattern of 5-7-5 is followed with the middle line longer than the other two.  The first line reflects $2r\pi$. In the design shown in Figure 6, the following relationship applies to chord $a$:   $a = 2r + 2m$ Looking at statement $q$ given in my *sangaku* formulae, the first relationship is $m$ related to $a$.  Based on $a = 2r + 2m$, the variable that $m$ shares with $a$ is $r$.

The geometry on the second line is  $2r + 2m + \frac{7}{9}\pi$.  Using the circle sector angle of 140°, line two contains the arc length $\left(\frac{7}{9}\pi\right)$ indicating the top of the second line's geometric figure; and the measure of the line at the bottom of the circle segment $(2r + 2m)$.  This expression clearly relates $r$ and $m$ as stated in the second relationship in the  $q$ formulae statement.  Notice the repetition of both 2 and $\pi$ in the first and second line of the haiku.

The third mathematical haiku line has $2s + 2s$.  Repeating $2s$ was chosen instead of $4s$ as it repeats the number 2, unifying all three haiku lines. The last portion of the $q$ formulae statement is the relationship of $s$ and $r$. As mentioned earlier, the first haiku geometric expression is in terms of $r$ with the last line focusing on $s$.

The square relates together the circle and the circle segment, thus creating a complete geometric thought; and with the geometric shapes, a complete geometric statement both visually and mathematically. All aspects of the Gion Shrine *sangaku* are unified; geometric design and geometric solution, and unified by the elegance found in simplicity.

The Gion Shrine *sangaku*'s answer was never supposed to be in numbers; but instead its inspiration rests in concepts. At the time when *wasanka* were discussing the need for more generalization, the Gion Shrine *sangaku* could only be solved using general geometric concepts. The Gion Shrine *sangaku* was designed to be different not just in its graceful, elegant design; but also different in its solution. This paper's solution fits the artistic nature of the Gion Shrine *sangaku*'s geometric image. The solution in geometry embeds the Gion Shrine *sangaku* deeply into the *sakoku* culture.

The Gion Shrine problem construction is simple but fastidiously precise, especially when viewed in terms of perimeter and circumference. It cannot be by chance that the Gion Shrine *sangaku* is also an unsolvable problem when traditional *sangaku* techniques are utilized. The author purposefully constructed the geometric figure so that it could not be solved by ordinary means; with intent, allowed for only one solution – a mathematical haiku - a solution in geometry.

Only one person painted the Gion Shrine *sangaku* design, inked the calligraphy and posed the problem.

I conclude that the author, although highly trained in mathematics when young, broke the tradition of his samurai family, and became one of the greatest masters of both haiku and *haiga* in all of Japanese history.  The Gion Shrine *sangaku* was not only about mathematics, it was about the Gion Shrine *sangaku* author's life and beliefs. I also believe that the Gion Shrine *sangaku* was a memorial to the author's mother on the twentieth anniversary of her death.

    Based on the Gion Shrine author's writings, I feel that there was also a broader statement being made about life in general. As mathematicians, we bring a synthesis of our life's experiences to every problem we approach. In the midst of solving a challenging problem, it is easy to forget to look outside of our training, and consider alternative solutions offered by the perspectives of other disciplines. Learning and living in cross-disciplines was an important enough concept to the Gion Shrine *sangaku* creator that he conveyed this sentiment to his disciples.  The more life an individual lives in every day, the more they bring to their craft.

# Chapter 5

## The Gion Shrine Problem Author

In 1716, in the small village of Kema, on the river Yodo, the Gion Shrine author was born under the surname Taniguchi. Little is known about Taniguchi's early life, although the surname Taniguchi would only have been given to a family of either samurai lineage or of significant prosperity.[65] Raised in wealth, and possibly a samurai family, he would have received significant mathematics education. Although actual documentation does not exist, experts believe that Taniguchi's mother died in 1728. This loss appears to have been enormously traumatic to Taniguchi. I surmise that this loss was part of his motivation for creating the Gion Shrine *sangaku*. Before the age of 21, Taniguchi left Kema, never to return.[66] Later in life, he changed his name to Yosa Buson. Buson is considered one of the greatest Japanese haiku poets that has ever lived. He was also considered a master artist. Perhaps the Gion Shrine problem was a reference to Buson's childhood, of which he never spoke.[67]

Constantly reading Chinese literature his entire life, Buson was well versed in all aspects of the language. As stated earlier, *sangaku* were scripted in *kanbun* a language that relates directly to Chinese characters and Chinese grammar. Later in life, when he decided to become a master artist, he studied Chinese painting style by reading a Chinese book of verse, *Mustard Seed Garden Manual of Painting*.[68]  Although his early poems reflected evidence of his Chinese studies, it became even more apparent as Buson aged. Towards the end of his life, his haiku not only showed the aesthetics of ancient Chinese poetry, but Buson dropped using Japanese *kana* script altogether in some haiku.[69]

---

[65] M. Ueda, *The Path of the Flowering Thorn: The Life and Poetry of Yosa Buson;* Stanford Univ. Press, CA 1998 p 3
[66] M. Ueda, *The Path of the Flowering Thorn: The Life and Poetry of Yosa Buson;* Stanford Univ. Press, CA 1998 p 2
[67] M. Ueda, *The Path of the Flowering Thorn: The Life and Poetry of Yosa Buson;* Stanford Univ. Press, CA 1998 p 112

Buson felt that haiku should be created to speak to all of the social classes. He also felt that there should be a synthesis among various fields of study.  "Buson's *rizoku* [rural custom] ideal was his resolution of a problem that had long been a site of debate in *haikai* [a broader term used for literature related to haiku]: how to balance the elements of high and low culture, the elite and the popular, that *haikai* brought together. '*Haikai* poets', he wrote, 'should read Chinese verse, painters should put down the brush and read books'…" [70]

"…a seeker of ideals that were more aesthetic than religious or moral…" is a description that Makoto Ueda gives Yosa Buson in *The Path of the Flowering Thorn: The Life and Poetry of Yosa Buson*.[71]  "…[Buson] freely let himself wander into the land of exotic beauty far removed from contemporary society and indulged in otherworldly dreams to his heart's content."

Certain that future aspiring haiku poets would read his writings, Buson recorded some conversations that he had with various disciples. In a conversation that Buson

---

[68] M. Ueda, *The Path of the Flowering Thorn: The Life and Poetry of Yosa Buson;* Stanford Univ. Press, CA  1998 p 67
[69] M. Ueda, *The Path of the Flowering Thorn: The Life and Poetry of Yosa Buson;* Stanford Univ. Press, CA  1998 p 17
[70] C.A. Crowley, H*aikai Poet Yosa Buson (1716-1783) and the Back to Basho Movement*; PhD. Thesis, Columbia University, 2001 p
[71] M. Ueda, *The Path of the Flowering Thorn: The Life and Poetry of Yosa Buson;* Stanford Univ. Press, CA  1998 p *v*

shared with one of his followers, Shōha, Buson wrote that he offered the following advice on writing beautiful haiku.[72]

"… Go far away from the marketplace, stroll among the trees in the garden, hold a banquet by a mountain stream, and enjoy a conversation over sake wine. … you will keep a leisurely frame of mind and enjoy the beauty of nature. Close your eyes and try composing a verse. When you come up with one, open your eyes.  … You will be … in a trance. Then the fragrance of blossoms will come drifting in the breeze, and the image of the moon will be seen floating on the water. That is the realm of *haikai* you are to enter."

Buson strongly believed that the more of life's experiences you brought to your art, no matter what that experience was, the more you would be able to give in perfecting the end results.

After years of travel, in 1751, Buson permanently settled in Kyoto, the city where the Gion Shrine is located.[73] It was at this time, that he began to use the name Yosa; also the name of the peninsula northwest of Kyoto, where it is believed his mother was born.

His loving humor was well known. It has been suggested that the following haiku might have been written by Buson as a humorous self-portrait.

   bolt upright
   such eyes fixed on the flowing clouds
   a frog

Figure 8 reflects Buson's loving humor in the *haiga*/haiku book he created depicting some of his fellow poets.

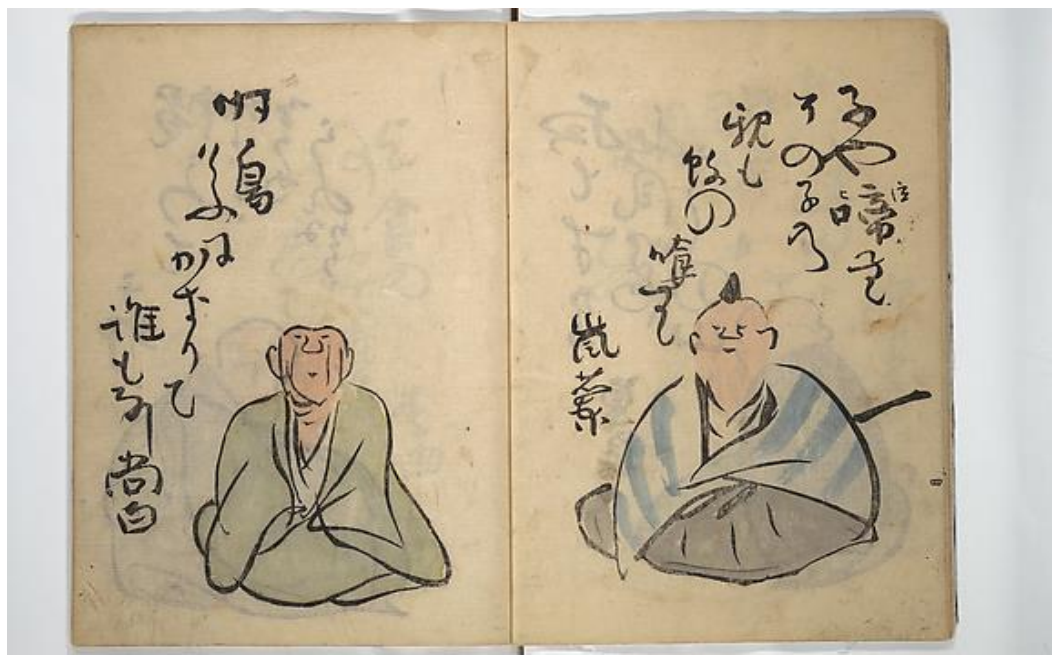
**Figure 8: "The Thirty-six Immortals of Haiku Verse" by Yosa Buson [27**

The question must be raised as to whether Buson created only one *sangaku*, the one that was hung at the Gion Shrine probably in 1748 (the twentieth anniversary of his mother's death). Could more *sangaku* be linked to him? With a dirth of information concerning Buson's early life; and a scant record of all *sangaku,* a recording compounded with many original tablets disappearing due to the passage of time, it is uncertain if this question will ever be answered. Other *sangaku* are simple by design, with some of them only having three geometric parts, but the *sangaku* on record do not seem to have the same grace and elegance as the one hung at the Gion Shrine in 1748. This question is left for others to answer.

# Chapter 6

# Justification

As mathematicians, surmising is not considered proof, especially when it is in connection to a mathematical problem that has gone unsolved for more than two and a half centuries. With historical documentation often being absent, and in this specific case, the original problem (the original *sangaku*) actually missing, there are some extrapolations that must be made. Many of the suppositions presented here are based on experts' interpretations and on as much documentation as was available close to when the Gion Shrine *sangaku* was written.

## 6.1    Justification of Solution

There are three pieces to the justification of the solution presented. The first rests in probable discussions occurring amongst Japanese mathematicians at the time that the Gion Shrine problem would have been inked on its wooden tablet. The second is focused on translation of ancient texts, and the variation of even current translations. Unification of image, the solution, and Japanese culture at the time is the third justification.

There is evidence that around the mid-eighteenth century, Japanese mathematicians began to discuss the need for *wasan* to move in a direction of generality. As mentioned in section 2.4, "*Wasan* - Japanese Mathematics", Ajima firmly recognized this need and made great strides in accomplishing that end. Although there are no direct records that Buson had direct contact with mathematicians during his life, he would have certainly been a part of the intellectual circles of the time, probably even as a samurai child. During the mid-1700s, he was traveling for portions of the time. It is very reasonable to assume that his travels would have included Kyoto as it was the major Japanese cultural center of the mid-1700s. Kyoto was also the Japanese center for mathematics at that time. It is very unlikely that he would not know about the *sangaku* at the various shrines and temples; or that he would not have heard the *wasanka* talking about the need for mathematical generality. The mathematical haiku solution is focused on very general aspects of geometry.

Translating a document written an hour ago from one language to another, can create misinterpretations simply due to the differences in languages. This is amplified when the translation involves the passage of hundreds of years. As noted earlier in this paper, comparing modern translations of the Gion Shrine *sangaku* indicates significant amounts of variation. The four modern interpretations range in date from 2008 to 2013 (a date could not be identified for the Cut-the-Knot web site post). With only a five year range in the interpretations, there is still variation due to different authors working on and researching the problem. As noted earlier, the use of the differing words is more than likely due to the use of either the 1815 rewriting of the original *sangaku* or the 1774 rewritten version. Even from that standpoint, there is only 41 years between the interpretation given by Mitsukuni (1815) and Ajima (1774). Language is a fluid, dynamical system. When mathematics heavily relies on language as a primary means of communication, it also becomes changeable and open to interpretation; or in this case, misinterpretation. I firmly believe that Buson, when he fashioned his beautiful *sangaku*, knew that misinterpretation would be likely. Buson had a lively sense of humor.

Rosalie Hosking in *Sangaku: A Mathematical, Artistic, Religious, and Diagrammatic Examination*[74] gives a thorough discussion of the opposing viewpoints when it comes to translating historical documents. The difficulty in accurate translation is even exaggerated when it comes to *sangaku* as they were written in *kanbun*, a noted cryptic language. Many of the current translations that accompany discussions of various *sangaku* are coached in current language, using terms that the *wasanka* would not have known and could not have used.[75] Hosking also amplifies the fact that *sangaku* were created with a definite link between the image present and the accompanying text. As pointed out by Hosking, this is particularly true when neither the problem description or answer have numerical values, as is the case with the Gion Shrine *sangaku*. Incorrectly translating the formulae on a *sangaku* leads to false conclusions. Partly due to a translation error, the entire meaning of the Gion Shrine *sangaku* was missed for hundreds of years.

*Sakoku* was a unique period of Japanese culture. Education of all types was encouraged allowing it to cut across social levels. While Japanese mathematics became specialized, other aspects of the Japanese culture found a rich environment for development. With the embracement of Confucianism, a shift occurred in the general

---

[74] R. Hosking, *Sangaku: A Mathematical, Artistic, Religious, and Diagrammatic Examination,* Univ. of Canterbury, UK, PhD. Thesis *https://ir.canterbury.ac.nz/handle/10092/12912* p 11-12

[75] R. Hosking, *Sangaku: A Mathematical, Artistic, Religious, and Diagrammatic Examination,* Univ. of Canterbury, UK, PhD. Thesis *https://ir.canterbury.ac.nz/handle/10092/12912* p 12

cultural personality towards directness. The simplicity of haiku folded perfectly into this cultural shift. Designed to be open-ended yet direct, haiku turned towards mathematics in the midst of *sakoku* fits profoundly well from a cultural standpoint. The concept of a mathematical haiku corresponds beautifully with the Edo Period culture.

The Gion Shrine *sangaku* was intended to be different not only in design but also in its solution; a solution that fits into every aspect of mid-eighteenth century Japanese culture.

## 6.2 Justification of the *Sangaku* Author

Justification for Yosa Buson as the Gion Shrine *sangaku* author is approached similar to a legal defense: opportunity, ability, and motive. Presenting these three aspects rests in further understanding Yosa Buson's studies, his life, and the words of advice that he gave to his haiku and *haiga* disciples.

Close to 1748, Buson was traveling. It is very reasonable to assume that his travels included Kyoto as it was the cultural center of Japan at that time. Records show that Buson made Kyoto his home in 1751.

Buson's constant study of Chinese material is well documented. Although other haiku poets might have also studied Chinese literature, Buson's level of passion for Chinese verse was not the norm. Partly for economic reasons, Buson decided to become recognized as an artist as well as a haiku poet. To accomplish this goal, he studied a Chinese book of verse, *Mustard Seed Garden Manual of Painting*, [76] instead of studying the works of his fellow Japanese artists. Buson later explained to his followers that he did this not to imitate Chinese painting, but to get himself in the correct frame of mind to create paintings. This reflects how greatly infused Chinese culture was with Buson. *Sangaku* are created in *kanbun* that has direct links to Chinese. Buson certainly would have had the expertise to intentionally construct a *kanbun* formulae for the Gion Shrine *sangaku* that would lead to multiple and differing interpretations. I contend that Buson was aware that this misinterpretation might occur. Unfortunately, the original tablet is no longer in existence. If it were still with us, it would display Buson's calligraphy.

---

[76] M. Ueda, *The Path of the Flowering Thorn: The Life and Poetry of Yosa Buson;* Stanford Univ. Press, CA 1998 p 67

Buson's genius is absolutely undeniable. His brilliance is unquestioned when it comes to both haiku and painting. He also might have had genius in mathematics.

Now consider the extreme precision of mind required to create a construction that has a bottom chord with the relationship of $a = 2r + 2m$. That same exact figure, when you locate the center point of the large circle, has the relationship of the large circle's radius equal to *2s*, twice the length of one side of the square. These precise relationships are <u>vital</u> to the mathematical haiku reflected by the Gion Shrine *sangaku*.

Mastering the art of creating haiku requires a precision of mental ability that is extreme. Haiku poets must drill down to the essence of ideas, images and situations. They must be able to construct and convey complete ideas in seventeen syllables.

Only the mind of an artist would visualize a distinct asymmetric design that had such perfect balance and grace as that of the Gion Shrine *sangaku* image. An artist would understand the need for related proportions in that same simple balanced design so that it would evoke a feeling of beauty.

The mental capacity to create a simple, elegant, beautiful problem containing meticulous design construction is the same mind required to strip naked aspects of life and lay down their essence in seventeen syllables. *Sangaku* and haiku are both means of communicating but through differing mechanisms.

Only a very rare individual could create a *sangaku* such as the one that hung at the Gion Shrine. Yosa Buson encapsulated all of the required characteristics.

Studying Buson's life, poetry, and his words of advice present motivation for his creation of the Gion Shrine *sangaku*.

Little is known about Yosa Buson's childhood from actual documentation or from his own words using direct references. There is record that his principle disciple, Takai Kitō (1741-1789) once described Buson's father as a "village headman". At a later time, this same manuscript was changed to give the description of Buson's father as "village gentleman". This description was then eliminated completely.[77] It is certain that Buson was the impetus behind Kitō's changes. Although much of his career, Buson went simply by the name Buson, he added "Yosa" in 1751. There was apparently some deep attachment to the use of this name. It "can be surmised that Buson's mother came from the Yosa peninsula northwest of Kyoto." [78] Buson experts acknowledge that Buson's attachment to his mother must have been great. In 1777, one year prior to the fiftieth anniversary of her death, Buson wrote several key nostalgic literary pieces that scholars connect to his childhood memories.

---

[77] M. Ueda, *The Path of the Flowering Thorn: The Life and Poetry of Yosa Buson;* Stanford Univ. Press, CA 1998 p 2
[78] M. Ueda, *The Path of the Flowering Thorn: The Life and Poetry of Yosa Buson;* Stanford Univ. Press, CA 1998 p 3

Included in the 1777 publication "Midnight Melodies", was an eighteen verse poem, *Kema Riverbank in the Spring Breeze*. This poem is written as if being told to Buson by a young woman. As Kema was his known birthplace, experts acknowledge that Buson was referring to his own feelings concerning his childhood. In a letter that Buson wrote in April, 1777, he refers to the poem as being "a journey scene from an amateurish play." [79] Following are some of the lines from *Kema Riverbank in the Spring Breeze* that are relevant to further discussions in this paper:

> "…
> as I walk down the bank to pick herbs
> thorn bushes block my way
> is it because they are jealous of me
> that they rip my clothes and scratch my thighs
> …" [80]

Later in that same year, 1777, *A New Florilegium* was published. In this volume were the following six haiku.

> "The Flower Festival –
> a mother's womb is only for
> temporary lodging"

> "eight of the fourth of the month –
> every baby died at birth
> is a Buddhist"

---

[79] M. Ueda, *The Path of the Flowering Thorn: The Life and Poetry of Yosa Buson;* Stanford Univ. Press, CA 1998 p 102

[80] M. Ueda, *The Path of the Flowering Thorn: The Life and Poetry of Yosa Buson;* Stanford Univ. Press, CA 1998 p 101-102

"time for summer clothes –

dew, shining white, starts to fall

on this body"

"time for summer clothes –

indeed, my mother had come

from the Fujiwara clan"

"a mountain cuckoo

the voice of a courtesan

composing a poem"

"hard of hearing

my father, a lay monk, cannot hear

the mountain cuckoo" [81]

"Some scholars believe that this sequence of six hokku [an older term used for haiku] has autobiographical implications. These poems …present images…making good sense only when it is tied to the poet's personal history." [82] Could Buson's mother have been used to produce children as Buson's father's legal wife was infertile? Maybe Buson's mother, talented in the arts, was treated like a concubine by Buson's father; a father who had no appreciation for the haiku that his mother created. It is possible that Buson's father with his "deafness" to the arts, was overzealous with Buson when it came to mathematics. Perhaps Buson had significant mathematical talent with a father

---

[81] M. Ueda, *The Path of the Flowering Thorn: The Life and Poetry of Yosa Buson;* Stanford Univ. Press, CA 1998 p 110

[82] M. Ueda, *The Path of the Flowering Thorn: The Life and Poetry of Yosa Buson;* Stanford Univ. Press, CA 1998 p 111

who was the jealous "thorn bush". These are suppositions, but similar to suppositions made by Buson scholars.

On the fiftieth anniversary of his mother's death, his nostalgia was documented through his own verse. It is certainly reasonable that he might have commemorated his mother's death at its twentieth anniversary, in 1748. Buson would have been close to 30 years old. The memory of his mother's passing would have still been in his mind, along with possible bitterness that he felt towards his father, and possibly towards mathematics itself.

By documenting the conversation that he had with his disciple Shōhu, Buson knew that the advice would be conveyed to all future haiku poets, and others such as ourselves. The message to be heard was to treat life as a banquet with every delectable relished with delight. He believed that each intellectual discipline gains from the others. Mathematics profits if *wasanka* try writing haiku or creating philosophy or charcoaling a drawing. Only then will you master your technique. Hanging a tablet with a mathematical problem solved only by thinking outside of mathematics itself, Buson crafted an image that evokes balance and grace; an image accompanied by a beautiful solution. Buson was sharing his art, his poetry, and his beliefs. Buson was a genius haiku poet, a master artist, and a brilliant mathematician.

Buson's message, '…stroll among the trees in the garden, hold a banquet by a mountain stream, and enjoy a conversation over sake wine. … the fragrance of blossoms will come drifting in the breeze, and the image of the moon will be seen floating on the water" [83] can still be heard today.

---

[83] M. Ueda, *The Path of the Flowering Thorn: The Life and Poetry of Yosa Buson;* Stanford Univ. Press, CA 1998 p 68

# Chapter 7

## Conclusion and
## *Sangaku* Author's Message

Mingling in Kyoto's cultural circles, Buson overhears several mathematicians discussing the need for mathematics to become more general in method and he embraces that thought. As the twentieth anniversary of his mother's passing approaches, he develops an idea that fits perfectly with his innate desire for the various intellectual disciplines to share their knowledge; for the disciplines' practitioners to learn from other aspects of intellectual pursuits. The *wasanka*, so focused on the existence of there being only one way to solve the *sangaku* moving in the breeze, forgot about the air itself. Awaiting discovery by the *wasanka,* there was another solution to the Gion Shrine *sangaku*; a solution that needed only a few numbers; a solution that rested in beauty; a solution that ends in simplicity.

Buson's message was one to be heard for generations. Einstein often suggested that everything should be challenged. With his thought experiments, he imagined the impossible. We all gained from his strange imaginings.

The past aides us well by setting a foundation to learn from in creating a new future.

It is important to find the simple forms that exist in complexity.

Imagine the impossible, and make it real.

*Gion Festival day*

*but there is a fragrant breeze*

*over Arrowroot Field* [84]

– Yosa Buson

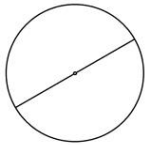

*the full harvest moon*

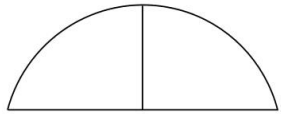

*just over the pond*

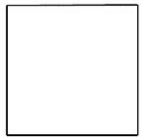

*filled by the rain* [85]

by Yosa Buson

---

[84] W.S. Merwin, T. Lento, *Collected Haiku of Yosa Buson;* Copper Canyon Press, Washington, 2013 p 113
[85] W.S. Merwin, T. Lento, *Collected Haiku of Yosa Buson;* Copper Canyon Press, Washington, 2013 p 143

**Images:**

**Web Sites:**